\documentclass[12pt]{amsart}
\usepackage{amssymb}

\usepackage{amsmath}

\theoremstyle{plain}

\numberwithin{equation}{section}
\input{tcilatex}

\begin{document}
\title[Geometric Equivalence]{On Geometrically Equivalent $S$-acts}
\author{Yefim Katsov}
\address{\textit{Department of Mathematics and Computer Science}\\
\textit{Hanover College, Hanover, IN 47243--0890, USA}}
\email{\textit{katsov@hanover.edu}}

\begin{abstract}
In this paper, considering the geometric equivalence for algebras of a
variety $_{S}\mathcal{A}$ of $S$-acts over a monoid $S$, we obtain
representation theorems describing all types of the equivalence classes of
geometrically equivalent $S$-acts of varieties $_{S}\mathcal{A}$ over groups 
$S$.
\end{abstract}

\maketitle

\section{Introduction}

{{I{n the last decade,} by Prof. B. I. Plotkin and his school, {algebraic
geometry} was introduced{\ in }a quite general setting --- namely, algebraic
geometry over algebras of an arbitrary variety of universal algebras (not
only of group varieties {as{,{{{{{{\ {for example, }}}}}}}in \cite%
{baummysrem:agog}}}) --- that brought forth a fascinating {new} area of
algebra known as universal algebraic geometry {{{{{(see,{{{{{{{\ } \textit{%
e.g.} }}}\cite{plotkin:voaaavcoav}, \cite{plotkin:snoagiua}, \cite%
{berzinsplts:agivoaw}, \cite{plotkin:slotuag}, and \cite%
{katlipplot:aocofmfsflm} for details}}})}}}}}. And }}one of the principal
problems in algebraic geometry over algebras of a variety $\Theta $ of
universal algebras involves studying interrelations between relations
between algebras $G_{1}$ and $G_{2}$ of $\Theta $ and relations between $%
G_{1}$- and $G_{2}$-geometries over them. As has repeatedly been emphasized
in \cite{plotkin:piaibuag},\cite{plotkin:slotuag},\cite{plotkin:snoagiua},
and \cite{plotkin:voaaavcoav}, the fundamental notion of \textit{geometric
equivalence}, heavily based on congruence theories of finitely generated
free algebras of the variety $\Theta $, proves to be crucial and important
in all investigations concerning this problem (one may consult \cite%
{plplts:geog}, \cite{berzins:geoa}, \cite{gobsh:rappcgeg}, \cite%
{tsurkov:geontfg}, and \cite{plts:atgeorog} for activity, obtained results
and open problems in this direction).

In this paper, continuing studying problems associated with algebraic
geometry over algebras of varieties $_{S}\mathcal{A}$ of $S$-acts over
monoids $S$, initiated in \cite{kat:apop}, we consider the geometric
equivalence of algebras in the varieties $_{S}\mathcal{A}$. After
introducing some notions and facts of algebraic geometry over algebras of $%
_{S}\mathcal{A}$ needed in a sequence, in Section 3, among other results, we
obtain the main results of the paper --- the representation theorems
(Theorems 3.21 and 3.22) --- describing all types of the equivalence classes
of geometrically equivalent $S$-acts of the varieties $_{S}\mathcal{A}$ over
groups $S$. We conclude the paper by stating three open problems,
delineating quite interesting and promising, in our view, directions for
further investigations.

Finally, all notions and facts of categorical algebra, used here without any
comments, can be found in \cite{macl:cwm}; for notions and facts from
universal algebra, we refer to \cite{gratzr:ua}.

\section{Basic preliminary notions and facts of algebraic geometry over $S$%
-acts}

\subsection{Algebraic varieties and closed congruences.}

Recall (see, \textit{e.g.}

\noindent \cite{klpknmx:maac}) that a \emph{left\/} $S$-\emph{act} over a
monoid $S$ is a non-empty set $A$ together with a scalar multiplication $%
(s,a)\mapsto sa$ from $S\times A$ to $A$ such that $1_{S}a=a$ and $%
(st)a=s(ta)$ for all $\ s,t\in S$ and $a\in A$. \emph{Right }$S$-\emph{act}
over $S$ and homomorphisms between $S$-acts are defined in the standard
manner. And, from now on, let $_{S}\mathcal{A}$ denote the category of left $%
S$-acts over a monoid $S$. Then (see, for example, \cite{klpknmx:maac}), a 
\textit{free} (left) $S$-act with a set of free generators, or a basis set, $%
I$ is a coproduct $\coprod_{i\in I}S_{i},S_{i}\cong $ $_{S}S$, $i\in I$, of
the copies of $_{S}S$ in the category $_{S}\mathcal{A}$. It is obvious that
two free $S$-acts $\coprod_{i\in I}S_{i}$ and $\coprod_{j\in J}S_{j}$ are
isomorphic in $_{S}\mathcal{A}$ iff $|I|=|J|$, and, hence, by %
\cite[Definition 2.8]{katlipplot:aocofmfsflm} $_{S}\mathcal{A}$ is an
IBN-variety. Let $X_{0}=\{x_{1},x_{2},\ldots ,x_{n},\ldots \}$ $\subseteq $ $%
\mathcal{U}$ be a fixed denumerable set of an infinite universe $\mathcal{U}$%
, and $_{S}\mathcal{A}^{0}$ the full subcategory of\ free $S$-acts $F_{X}%
\overset{def}{=}\coprod_{x\in X}S_{x},$ $X\ \subseteq \mathcal{U},|X|<$ $%
\infty $, with a finite basis $X$ of the variety $_{S}\mathcal{A}$.

For a fixed $S$-act $G\in |_{S}\mathcal{A}|$ of the variety $_{S}\mathcal{A}$
and a free algebra $F_{X}\in |_{S}\mathcal{A}^{0}|$, the set $_{S}\mathcal{A}%
(F_{X},G)$ of homomorphisms from $F_{X}$ to $G$ is treated as an affine
space over $G$ consisting of points (homomorphisms) $\mu
:F_{X}\longrightarrow G$. Any system $T=\{s_{x_{i}}=t_{y_{i}}|$ $%
x_{i},y_{i}\in X,s_{x_{i}}\in $ $S_{x_{i}},t_{y_{i}}\in $ $S_{y_{i}},i\in
I\} $ of the equations in free variables of $X$ in $F_{X}$ can be obviously
viewed as the binary relation $T\subseteq F_{X}\times F_{X}$ on $F_{X}$, and
a point (homomorphism) $\mu :F_{X}\longrightarrow G$ is a \textit{solution}
of $T$ iff $\mu (s_{x_{i}})=\mu (t_{y_{i}})$ for any equation $%
s_{x_{i}}=t_{y_{i}}$ of $T$, \textit{i.e.}, iff $T\subseteq Ker$ $\mu $.
Thus, for any set of points $A\subseteq $ $_{S}\mathcal{A}(F_{X},G)$ and any
binary relation $T\subseteq F_{X}\times F_{X}$ on $F_{X}$, the assignments%
\begin{equation*}
T\longmapsto T_{G}^{^{\prime }}\overset{def}{=}\{\mu :F_{X}\longrightarrow
G|T\subseteq Ker\mu \}\text{ and }A\longmapsto A^{^{\prime }}\overset{def}{=}%
\cap _{\mu \in A}Ker\mu
\end{equation*}%
define the Galois correspondence between binary relations (or systems of
equations) $T\ $on $F_{X}$ and sets of points $A\ $of the space $_{S}%
\mathcal{A}(F_{X},G)$. A congruence $T$ on $F_{X}$ is said to be $G$-\textit{%
closed} if $T=A^{^{\prime }}$ for some point set $A\subseteq $ $_{S}\mathcal{%
A}(F_{X},G)$; a point set $A$ is called a $G$-\textit{closed set}, or\textit{%
\ }an \textit{algebraic variety}, in the space $_{S}\mathcal{A}(F_{X},G)$ if 
$A=T_{G}^{^{\prime }}$ for some relation $T\ $on $F_{X}$. As usual, the
Galois correspondence produces the closures: $A^{^{\prime \prime }}\overset{%
def}{=}(A^{^{\prime }})^{^{\prime }}$and $T^{^{\prime \prime }}\overset{def}{%
=}(T_{G}^{^{\prime }})^{^{\prime }}$. And one easily obtains the following
observation.\medskip

\noindent \textbf{Proposition 2.1.} \textit{A congruence} $T\subseteq
F_{X}\times F_{X}$ \textit{on} $F_{X}$\textbf{\ }\textit{is }$G$\textit{%
-closed }

\noindent \textit{iff }$T=T^{^{\prime \prime }}$\textit{.}\textbf{\ \medskip 
}

\noindent \textbf{Proof}. $\Longrightarrow $. Let $T=A^{^{\prime }}$ for a
point set $A\subseteq $ $_{S}\mathcal{A}(F_{X},G)$, \textit{i.e.}, $T=\cap
_{\mu \in A}Ker\mu $. Then, $T_{G}^{^{\prime }}=\{\lambda
:F_{X}\longrightarrow G|\cap _{\mu \in A}Ker\mu =T\subseteq Ker\lambda \}$,
and hence, $T^{^{\prime \prime }}=\cap _{\lambda \in T_{G}^{^{\prime
}}}Ker\lambda =\cap _{\mu \in A}Ker\mu \cap _{\lambda \in T_{G}^{^{\prime
}}\diagdown A}Ker\lambda =\cap _{\mu \in A}Ker\mu =T$.

$\Longleftarrow $. $T=T^{^{\prime \prime }}=(T_{G}^{^{\prime }})^{^{\prime
}}=A^{^{\prime }}$, where $A=$ $T_{G}^{^{\prime }}$.\textit{\ \ \ \ \ \ }$%
_{\square }\medskip $

\noindent \textbf{Corollary 2.2.} (\textit{cf. }\cite[Proposition 2.1]%
{plotkin:slotuag}) \textit{A congruence} $T\subseteq F_{X}\times F_{X}$ 
\textit{on} $F_{X}$\textbf{\ }\textit{is }$G$\textit{-closed iff \ there is
an embedding }$\mu :$\textit{\ }$F_{X}/T\rightarrowtail G^{I}$ \textit{for
some set }$I$\textit{.}\textbf{\ \medskip }

\noindent \textbf{Proof}. $\Longrightarrow $. Since $T=T^{^{\prime \prime
}}=(T_{G}^{^{\prime }})^{^{\prime }}$, by \cite[Theorem 20.2]{gratzr:ua} $%
F_{X}/T$ is a subdirect product of the $S$-subacts $F_{X}/Ker\lambda
\subseteq G$ of $_{S}G$, $\lambda \in T_{G}^{^{\prime }}$, and therefore,
there exists an embedding $F_{X}/T\rightarrowtail G^{I}$ for the set $I=$%
\textit{\ }$T_{G}^{^{\prime }}$.

\textit{\ }$\Longleftarrow $. If \textit{\ }$\mu :$\textit{\ }$%
F_{X}/T\rightarrowtail G^{I}$, and $\pi _{i}:G^{I}\twoheadrightarrow G$, $%
i\in I$, are the canonical projections, then $T=\cap _{i\in I}Ker$ $\pi
_{i}\mu $.\textit{\ \ \ \ \ \ }$_{\square }\medskip $

Analogously, we have\medskip

\noindent \textbf{Proposition 2.3.} \textit{A point set }$A\subseteq $ $_{S}%
\mathcal{A}(F_{X},G)$ \textit{is an algebraic variety }

\noindent \textit{iff} $A=A^{^{\prime \prime }}$\textit{.}\textbf{\ \medskip 
}

\noindent \textbf{Proof}. $\Longrightarrow $. Let $A=T_{G}^{^{\prime
}}=\{\mu :F_{X}\longrightarrow G|$ $T\subseteq Ker\mu \}$ for some relation $%
T\subseteq F_{X}\times F_{X}$ on $F_{X}$. Then, $A^{^{\prime }}=\cap _{\mu
\in A}Ker\mu \supseteq T$ and, hence, $A^{^{\prime \prime }}=\{\lambda
:F_{X}\longrightarrow G|$ $T\subseteq \cap _{\mu \in A}Ker\mu \subseteq
Ker\lambda \}$, and therefore, $A^{^{\prime \prime }}\subseteq A\subseteq
A^{^{\prime \prime }}$.

$\Longleftarrow $. Obvious, as $A=A^{^{\prime \prime }}=(A^{^{\prime
}}=T)^{^{\prime }}$.\textit{\ \ \ \ \ \ }$_{\square }\medskip $

\subsection{The functors $Alv_{G}$ and $Cl_{G}$ over $_{S}\mathcal{A}^{0}$.}

We consider two natural functors: $Alv_{G}:$ $_{S}\mathcal{A}%
^{0}\longrightarrow Set$ and $Cl_{G}:$ $(_{S}\mathcal{A}^{0})^{op}%
\longrightarrow Set$.

a) $Alv_{G}:$ $_{S}\mathcal{A}^{0}\longrightarrow Set$. For any $F_{X}\in
|_{S}\mathcal{A}^{0}|$, $Alv_{G}(F_{X})\overset{def}{=}\{$ $A|$ $A\subseteq $
$_{S}\mathcal{A}(F_{X},G)$ $\&$ $A=A^{^{\prime \prime }}\}$; and for any
homomorphism $s\in $ $_{S}\mathcal{A}^{0}(F_{Y},F_{X})$ and $B=\{\mu
:F_{Y}\longrightarrow G|$ $B^{^{\prime }}\subseteq Ker\mu \}\in $ $%
Alv_{G}(F_{Y})$, $\ Alv_{G}(s)\overset{def}{=}s^{-1}B\overset{def}{=}%
\{\alpha :F_{X}\longrightarrow G|$ $\alpha s\in B\}=\{\alpha
:F_{X}\longrightarrow G|$ $sB^{^{\prime }}\subseteq Ker\alpha \}$, where the
relation $sB^{^{\prime }}\overset{def}{=}\{$ $(\omega ,\omega ^{^{\prime
}})\in F_{X}\times F_{X}|$ $\exists (\omega _{0},\omega _{0}^{^{\prime
}})\in B^{^{\prime }}\subseteq F_{Y}\times F_{Y}:s(\omega _{0})=\omega $ $\&$
$s(\omega _{0}^{^{\prime }})=\omega ^{^{\prime }}\}$.\smallskip

b) $Cl_{G}:$ $(_{S}\mathcal{A}^{0})^{op}\longrightarrow Set$. For any $%
F_{X}\in |_{S}\mathcal{A}^{0}|$, $Cl_{G}(F_{X})\overset{def}{=}\{$ $T|$ $%
T\subseteq F_{X}\times F_{X}$ $\&$ $T=T^{^{\prime \prime }}\}$; and for any
homomorphism $s\in $ $_{S}\mathcal{A}^{0}(F_{Y},F_{X})$ and $T=T^{^{\prime
\prime }}=(T_{G}^{^{\prime }})^{^{\prime }}=\cap _{\alpha \in
T_{G}^{^{\prime }}}Ker\alpha \subseteq F_{X}\times F_{X}$, $\ Cl_{G}(s)%
\overset{def}{=}s^{-1}T\overset{def}{=}\cap _{\alpha \in T_{G}^{^{\prime
}}}Ker$ $s\alpha \subseteq F_{Y}\times F_{Y}$.\smallskip

One can easily see that the just defined mappings $Alv_{G}$ and $Cl_{G}$ are
indeed covariant and contravariant functors, respectively. Also, from the
results of \cite[Section 2.2]{plotkin:slotuag} it is easy to see that for
any $F_{X}\in |_{S}\mathcal{A}^{0}|$ the sets $Alv_{G}(F_{X})$ and $%
Cl_{G}(F_{X})$ are meet-semilattices with the respect to the natural partial
orders on the set of all subsets of $_{S}\mathcal{A}(F_{X},G)$ and $%
F_{X}\times F_{X}$, respectively. Moreover, if we define $A\overline{\cup }B%
\overset{def}{=}(A\cup B)^{^{\prime \prime }}$and $T_{1}\overline{\cup }T_{2}%
\overset{def}{=}(T_{1}\cup T_{2})^{^{\prime \prime }}$for any $A,B\in $ $\
Alv_{G}(F_{X})$ and $T_{1},T_{2}\in Cl_{G}(F_{X})$, the sets $Alv_{G}(F_{X})$
and $Cl_{G}(F_{X})$ become lattices, and by \cite[Proposition 2.4]%
{plotkin:slotuag} (see also \cite[Proposition 4]{plotkin:snoagiua}) the
assignments $T\longmapsto T^{^{\prime }}$ for each $T\in Cl_{G}(F_{X})$
establish a dual isomorphism between the lattices $Cl_{G}(F_{X})$ and $%
Alv_{G}(F_{X})$. Following \cite[Definition 2.4]{plotkin:slotuag}, we call
an $S$-act $G$ \textit{geometrically stable} iff $A\overline{\cup }B=A\cup B$
for any $F_{X}\in |_{S}\mathcal{A}^{0}|$ and $A,B\in $ $\ Alv_{G}(F_{X})$,
or, as one can easily see by using the results of \cite[Section 2.2]%
{plotkin:slotuag}, iff $(T_{1}^{^{\prime }}\cup T_{2}^{^{\prime
}})=(T_{1}\cap T_{2})^{^{\prime }}$for any $F_{X}\in |_{S}\mathcal{A}^{0}|$
and $T_{1},T_{2}\in Cl_{G}(F_{X})$. Thus, considering a particular case of
the varieties $_{S}\mathcal{A}$, namely the case when the monoid $S=\{1\}$,
and therefore, $_{S}\mathcal{A}$ is just the category $Set$ of non-empty
sets, and in contrast to \cite[Theorem 1]{plotkin:snoagiua}, we have the
following observation.\medskip

\noindent \textbf{Proposition 2.4.} \textit{In }$Set$\textit{, the
singletons are the only geometrically stable sets.}\textbf{\ \medskip }

\noindent \textbf{Proof}. Let $G,X\in Set$, and $|G|=1$. As it is clear that
all congruences on $X$ are just equivalence relations, $Cl_{G}(X)$ has only
the universal equivalence relation, and $|Alv_{G}(X)|=1$, $G$ is
geometrically stable.

Let $G,X\in Set$, $|G|\geq 2$, and $T$ an equivalence relation on $X$. Then
from obvious and well-known set theory facts, there are embeddings $%
X/T\rightarrowtail 2^{|X/T|}\rightarrowtail G^{^{|X/T|}}$, and therefore, by
Corollary 2.2 the equivalence relation $T$ is $G$-closed. Now take $%
X=\{x_{1},x_{2},x_{3}\}$, and let $T_{1}=%
\{(x_{1},x_{1}),(x_{2},x_{2}),(x_{3},x_{3}),(x_{1},x_{2}),(x_{2},x_{1})\}$
and $T_{2}=\{(x_{1},x_{1}),(x_{2},x_{2}),$

\noindent $(x_{3},x_{3}),(x_{2},x_{3}),(x_{3},x_{2})\}$ be equivalence
relations on $X$. Then, $(T_{1}\cap
T_{2})=\{(x_{1},x_{1}),(x_{2},x_{2}),(x_{3},x_{3})\}$ and $(T_{1}\cap
T_{2})^{^{\prime }}=Set(X,G)$; and $T_{1}^{^{\prime }}=\{\alpha
:X\longrightarrow G$ $|$ $\alpha (x_{1})=\alpha (x_{2})\}$, $T_{2}^{^{\prime
}}=\{\beta :X\longrightarrow G$ $|$ $\beta (x_{2})=\beta (x_{3})\}$. For $%
g_{1},g_{2}\in G$, $g_{1}\neq g_{2}$, there exists $\gamma :X\longrightarrow
G$ such that $\gamma (x_{1})=\gamma (x_{3})=g_{1}$ and $\gamma (x_{2})=g_{2}$%
. Thus, $\gamma \in (T_{1}\cap T_{2})^{^{\prime }}$, but $\gamma \notin
(T_{1}^{^{\prime }}\cup T_{2}^{^{\prime }})$.\textit{\ \ \ \ \ \ }$_{\square
}\medskip $

\section{Geometrically equivalent $S$-acts}

For varieties of $S$-acts, the notion of \textit{geometric equivalence} is
defined as follows.\medskip

\noindent \textbf{Definition 3.1.} (\textit{cf.} \cite[Definition 3.1]%
{plotkin:slotuag} and \cite{plotkin:snoagiua}) Two $S$-acts $G_{1}$, $%
G_{2}\in |_{S}\mathcal{A}|\ $are said to be \textit{geometrically equivalent}%
,\textit{\ }$G_{1}\overset{\triangle }{\sim }G_{2}$, iff $%
T_{G_{1}}^{^{\prime \prime }}=T_{G_{2}}^{^{\prime \prime }}$ for any binary
relation $T\subseteq F_{X}\times F_{X}$ on any free $S$-acts $F_{X}\in |_{S}%
\mathcal{A}^{0}|$, \textit{i.e.}, iff $Cl_{G_{1}}(F_{X})=Cl_{G_{2}}(F_{X})$
for any free $S$-acts $F_{X}\in |_{S}\mathcal{A}^{0}|$, or iff $%
Cl_{G_{1}}=Cl_{G_{2}}\ $as functors.\medskip\ \ \ \ \ \ \ \ \ \ \ \ \ \ \ \
\ \ \ \ \ \ \ \ \ \ \ \ \ \ \ \ \ \ \ \ \ \ \ \ \ \ \ \ \ \ \ \ \ \ \ \ \ \
\ \ \ \ \ \ \ \ \ \ \ \ \ \ \ \ \ \ \ \ \ \ \ \ \ 

From the proof of Proposition 2.4 one immediately obtains the following
description of geometrically equivalent classes in $Set$.\medskip

\noindent \textbf{Theorem 3.2.} \textit{In} $Set$\textit{, there are only
two classes of geometrically equivalent sets: singletons, and all sets with
more then one element.\ \ \ \ \ \ }$_{\square }\medskip $

\noindent \textbf{Lemma 3.3.} \textit{For any relation} $T\subseteq
F_{X}\times F_{X}$ \textit{and monomorphism} $\mu :G_{1}\rightarrowtail
G_{2} $, $T_{G_{1}}^{^{\prime \prime }}\supseteq T_{G_{2}}^{^{\prime \prime
}}$\textit{.}\textbf{\ \medskip }

\noindent \textbf{Proof}. Indeed, $T_{G_{1}}^{^{\prime }}=\{\alpha
:F_{X}\longrightarrow G_{1}$ $|$ $T\subseteq Ker$ $\alpha =Ker$ $\mu \alpha
\}$, $T_{G_{2}}^{^{\prime }}=\{\beta :F_{X}\longrightarrow G_{2}$ $|$ $%
T\subseteq Ker\beta \}\supseteq \{\mu \alpha :F_{X}\longrightarrow G_{2}$ $|$
$\alpha :F_{X}\longrightarrow G_{1}$ $\&$ $T\subseteq Ker$ $\alpha \}$.
Therefore, $T_{G_{1}}^{^{\prime \prime }}=\cap _{\alpha \in
T_{G_{1}}^{^{\prime }}}Ker$ $\alpha =\cap _{\alpha \in T_{G_{1}}^{^{\prime
}}}Ker$ $\mu \alpha \supseteq \cap _{\beta \in T_{G_{2}}^{^{\prime }}}Ker$ $%
\beta =T_{G_{2}}^{^{\prime \prime }}$.\textit{\ \ \ \ \ \ }$_{\square }$

\noindent \textbf{Corollary 3.4.} (\cite[Proposition 3.1]{plotkin:slotuag}) $%
G\overset{\triangle }{\sim }G^{I}$ \textit{for any} $G$ \textit{and any set }%
$I$\textit{.}\textbf{\ \medskip }

\noindent \textbf{Proof}. Because of the diagonal embedding $%
G\rightarrowtail G^{I}$, it follows from Lemma 3.3 that $T_{G}^{^{\prime
\prime }}\supseteq T_{G^{I}}^{^{\prime \prime }}$ for any relation $%
T\subseteq F_{X}\times F_{X}$.

Now let $\pi _{i}:G^{I}\twoheadrightarrow G_{i}$, $G_{i}=G,$ $i\in I,$ be
the natural projections, $T$ a binary relation on\ $F_{X}$, and $\alpha \in $
$T_{G^{I}}^{^{\prime }}$. Then, $T\subseteq Ker$ $\pi _{i}\alpha $ for any $%
i\in I$, and, hence, $\pi _{i}\alpha \in $ $T_{G}^{^{\prime }},$ $i\in I$.
And it is clear that $Ker$ $\alpha =\cap _{i\in I}Ker$ $\pi _{i}\alpha
\supseteq \cap _{\beta \in T_{G}^{^{\prime }}}Ker$ $\beta =T_{G}^{^{\prime
\prime }}$, and therefore, $T_{G^{I}}^{^{\prime \prime }}=\cap _{\alpha \in
T_{G^{I}}^{^{\prime }}}Ker$ $\alpha \supseteq T_{G}^{^{\prime \prime }}$.%
\textit{\ \ \ \ \ \ }$_{\square }$

The next observation will prove to be a quite useful generalization of
Corollary 3.4.\medskip\ 

\noindent \textbf{Corollary 3.5.} \textit{Let for }$G_{1}$ \textit{and }$%
G_{2}$\textit{\ there exist sets} $I$ \textit{and} $J$ \textit{and embeddings%
} $G_{1}\rightarrowtail G_{2}^{I}$ \textit{and} $G_{2}\rightarrowtail
G_{1}^{J}$\textit{. Then,} $G_{1}\overset{\triangle }{\sim }G_{2}$\textit{.\
\ \ \ \ \ }$_{\square }\medskip $

Now let $_{S}\mathcal{B}\subseteq $ $_{S}\mathcal{A}$ be a subvariety of the
variety $_{S}\mathcal{A}$. As is well-known (see, for example, \cite[Chapter
VI]{mal:as}), a free algebra with a finite basis $X$ in the variety $_{S}%
\mathcal{B}$ is, in fact, the factor algebra $F_{X}/\tau $ of $F_{X}\in |_{S}%
\mathcal{A}^{0}|$ modulo $\tau $, where $\tau $ is a corresponding verbal
congruence. For a congruence $T$ on $F_{X}/\tau $ and the epimorphism $\pi
:F_{X}\twoheadrightarrow F_{X}/\tau $ corresponding $\tau $, let $\pi ^{-1}T%
\overset{def}{=}\{(\omega ,\omega ^{^{\prime }})\in F_{X}\times F_{X}|$ $%
(\pi (\omega ),\pi (\omega ^{^{\prime }}))\in T\}$ be the inverse image of
the congruence $T$. In these notations, we have the following
observation.\medskip

\noindent \textbf{Proposition 3.6.} \textit{For any} $G\in |_{S}\mathcal{B}|$%
\textit{, the mapping} $T\longmapsto \pi ^{-1}T$ \textit{is a bijection
between the sets of} $G$\textit{-closed congruences in} $_{S}\mathcal{B}$ 
\textit{and in} $_{S}\mathcal{A}$\textit{.} \textbf{\ \medskip }

\noindent \textbf{Proof}. Let $T$ be $G$-closed in $_{S}\mathcal{B}$. Then, $%
T_{G}^{^{\prime }}=\{\alpha :F_{X}/\tau \longrightarrow G$ $|$ $T\subseteq
Ker$ $\alpha \}$ and $T=T_{G}^{^{\prime \prime }}=\cap _{\alpha \in
T_{G}^{^{\prime }}}Ker$ $\alpha $. One may easily see that the latter is a
limit of a corresponding diagram, and therefore, applying the results of %
\cite[Section 9.8]{macl:cwm} on iterated limits, obtain that $\pi
^{-1}T=\cap _{\alpha \in T_{G}^{^{\prime }}}Ker$ $\alpha \pi $, and, hence, $%
\pi ^{-1}T$ is $G$-closed in $_{S}\mathcal{A}$, too.

Now, if a congruence $T$ on $F_{X}$ is $G$-closed, then, since $G\in |_{S}%
\mathcal{B}|$, one has $T=$ $\cap _{\alpha \in A}Ker$ $\alpha \pi $ for some 
$A\subseteq $ $_{S}\mathcal{B}(F_{X}/\tau ,G)$, and, again taking into
consideration the iterated limits, obtains that $T=$ $\pi ^{-1}(\cap
_{\alpha \in A}Ker$ $\alpha )$. As $\cap _{\alpha \in A}Ker$ $\alpha $ is a $%
G$-closed congruence on $F_{X}/\tau $, we conclude that any $G$-closed
congruence on $F_{X}$ is an inverse image of $G$-closed congruence on $%
F_{X}/\tau $.

Noting that $\pi (\pi ^{-1}T)=T$ for any congruence $T$ on $F_{X}/\tau $, we
end the proof.\textit{\ \ \ \ \ \ }$_{\square }\medskip $

\noindent \textbf{Corollary 3.7.} (\cite[Proposition 3.2]{plotkin:slotuag}) 
\textit{For any two algebras} $G_{1},G_{2}\in |_{S}\mathcal{B}|$ \textit{of
a subvariety} $_{S}\mathcal{B}$ \textit{of the variety} $_{S}\mathcal{A}$, $%
G_{1}\overset{\triangle }{\sim }G_{2}$ \textit{in} $_{S}\mathcal{A}$ \textit{%
iff} $G_{1}\overset{\triangle }{\sim }G_{2}$ \textit{in} $_{S}\mathcal{B}$%
\textit{.\ \ \ \ \ \ }$_{\square }\medskip $

Defining an $S$-act $A$ to be \textit{trivial} iff $sa=a$ for any $s\in S$
and $a\in A$, one readily sees that all trivial $S$-acts form a subvariety
of the variety $_{S}\mathcal{A}$, defined by the identities $\forall
x(sx=x), $ $s\in S$, and actually coinciding with $Set$. From this
observation, Corollary 3.7, and Theorem 3.2\ we have\medskip

\noindent \textbf{Corollary 3.8.} \textit{In} $_{S}\mathcal{A}$\textit{, all
trivial }$S$\textit{-acts of cardinality greater than one are geometrically
equivalent.\ \ \ \ \ \ }$_{\square }\medskip $

As usual, an $S$-act $A\in |_{S}\mathcal{A}|$ is \textit{cyclic} if there is 
$a\in A$ such that $_{S}A=\{Sa\}$, and $_{S}A$ is called \textit{simple} if
it has no proper subacts. It is obvious that any simple act is cyclic.

From now on, we assume that a monoid $S$ is a group. Then, it is clear that
any cyclic act is simple, and, by \cite[Proposition 1.5.34]{klpknmx:maac},
any act $A\in |_{S}\mathcal{A}|$ is a disjoint union, or a coproduct in the
category $_{S}\mathcal{A}$, of simple (cyclic) subacts. Moreover, using %
\cite[Proposition 1.5.17]{klpknmx:maac}, one can easily see that for any
cyclic act $_{S}\overline{G}$ there exists a subgroup $G$ of $S$ such that $%
_{S}\overline{G}$ is isomorphic to the act $_{S}S/G$ of the left cosets
modulo $G$; and, therefore, any\ $A\in |_{S}\mathcal{A}|$, in fact, can be
considered as a coproduct of suitable left cosets $_{S}S/G\overset{def}{=}$ $%
_{S}\overline{G}$ of the group $S$.\medskip

\noindent \textbf{Proposition 3.9.} \textit{Let }$G$\textit{\ and }$H$%
\textit{\ be subgroups of a group }$S$\textit{, and the cyclic }$S$\textit{%
-acts }$_{S}S/G=$ $_{S}\overline{G}$\textit{\ and }$_{S}S/H=$ $_{S}\overline{%
H}$\textit{\ geometrically equivalent. Then there are }$\alpha ,\beta \in S$%
\textit{\ such that }$G^{\alpha }\overset{def}{=}\alpha ^{-1}G\alpha
\subseteq H$\textit{\ and }$H^{\beta }\overset{def}{=}\beta ^{-1}H\beta
\subseteq G$\textit{.} \textbf{\ \medskip }

\noindent \textbf{Proof}. Let $T_{1}=\{(sg_{1},sg_{2})$ $|$ $g_{1,}g_{2}\in
G $, $s\in S\}$, $T_{2}=\{(sh_{1},sh_{2})$ $|$ $h_{1,}h_{2}\in H$, $s\in S\}$
be the congruences on the act $_{S}S$, and $\pi _{1}:$ $_{S}S\longrightarrow 
$ $_{S}S/T_{1}=$ $_{S}S/G=$ $_{S}\overline{G}$, $\pi _{2}:$ $%
_{S}S\longrightarrow $ $_{S}S/T_{2}=$ $_{S}S/H=$ $_{S}\overline{H}$ the
canonical epimorphisms corresponding to the subgroups $G,H\subseteq S$,
respectively.

As $T_{1}=Ker$ $\pi _{1}$, $T_{1}\in $\ $Cl_{\overline{G}}(_{S}S)$, and,
since $\overline{G}\overset{\triangle }{\sim }\overline{H}$, $T_{1}\in $\ $%
Cl_{\overline{H}}(_{S}S)$, too. Therefore, there exists a homomorphism $f:$ $%
_{S}S\longrightarrow $ $_{S}\overline{H}$ such that $Ker$ $f\supseteq $ $%
T_{1}$, and, hence, there is such a homomorphism $\overline{\alpha }:$ $_{S}%
\overline{G}\longrightarrow $ $_{S}\overline{H}$ that $f=$ $\overline{\alpha 
}\pi _{1}$. From this and noting that the free object $_{S}S\in |_{S}%
\mathcal{A}|$ is obviously projective in the category $_{S}\mathcal{A}$
(see, also \cite[Proposition III.17.2]{klpknmx:maac}), we have $\overline{%
\alpha }\pi _{1}=\pi _{2}\alpha $ for some $\alpha :$ $_{S}S\longrightarrow $
$_{S}S$. Agreeing to write endomorphisms of $_{S}S$ on the right of the
elements they act on, one can easily see that actions of endomorphisms
actually coincide with multiplications of elements of $_{S}S$ on the right
by elements of the monoid $S$, and, therefore, $End(_{S}S)=$ $S$, and, thus, 
$\forall s\in S:\alpha (s)\overset{def}{=}s\alpha $, where $\alpha \in S$.\
Then, from the equation $\overline{\alpha }\pi _{1}=\pi _{2}\alpha $, we
have: $(sg_{1},sg_{2})\in $ $T_{1}\Longrightarrow (sg_{1}\alpha
,sg_{2}\alpha )\in $ $T_{2}$ $\Longleftrightarrow $ $g_{1}\alpha
h=g_{2}\alpha $ for some $h\in H$ $\Longleftrightarrow $ $(\alpha
^{-1}g_{1}\alpha )h=\alpha ^{-1}g_{2}\alpha $ for some $h\in H$. Therefore,
we have that $G^{\alpha }=\alpha ^{-1}G\alpha \subseteq tH$ for some $t\in S$%
. Show that, in fact, $t\in H$. Indeed, if $g_{1,}g_{2}\in G$, then there
exist $h_{1,}h_{2},h_{3}\in H$ such that $\alpha ^{-1}g_{1}\alpha
=th_{1},\alpha ^{-1}g_{2}\alpha =th_{2},\alpha ^{-1}g_{1}g_{2}\alpha =th_{3}$%
, and, hence, $th_{3}=th_{1}th_{2}$ from which immediately follows that $%
t\in H$.\ Thus, we have proved an inclusion $G^{\alpha }=\alpha ^{-1}G\alpha
\subseteq H$ for some $\alpha \in S$.

By symmetry, we also obtain an inclusion $H^{\beta }=\beta ^{-1}H\beta
\subseteq G$ for some $\beta \in S$.\textit{\ \ \ \ \ \ }$_{\square
}\medskip \smallskip $

\noindent \textbf{Proposition 3.10.} \textit{Let }$G$\textit{\ and }$H$%
\textit{\ be subgroups of a group }$S$\textit{, and }$H\subset G$\textit{.
Then\medskip , }$_{S}\overline{H}\overset{\triangle }{\nsim }$\textit{\ }$%
_{S}\overline{G}$\textit{. }\textbf{\ }

\noindent \textbf{Proof}. Let $\pi _{\overline{G}}:$ $_{S}S\longrightarrow $ 
$_{S}\overline{G}$, $\pi _{\overline{H}}:$ $_{S}S\longrightarrow $ $_{S}%
\overline{H}$ be the canonical epimorphisms corresponding to the subgroups $%
G,H\subseteq S$, respectively. Then, $T=Ker$ $\pi _{\overline{H}}=\{(s,t)$ $%
| $ $s,t\in S$ $\&$ $sH=tH\}$ obviously is an $_{S}\overline{H}$-closed
congruence, \textit{i.e.}, $T=T_{\overline{H}}^{^{\prime \prime }}$.

Hence, assuming $_{S}\overline{H}\overset{\triangle }{\sim }$ $_{S}\overline{%
G}$, we\ have that $T$ is an $_{S}\overline{G}$-closed congruence, too; and
therefore, there exist a family of homomorphisms $f_{i}:$ $%
_{S}S\longrightarrow $ $_{S}\overline{G},i\in I,$ such that $T=\cap _{i\in
I}Ker$ $f_{i}$. However, using the Yoneda lemma (see, for example, %
\cite[Section III.2]{macl:cwm}) and the fact that $_{S}S$ is a
projective(free) object in the category $_{S}\mathcal{A}$, one immediately
has that for any $i\in I$, there exists an isomorphism $\alpha _{i}\in
Iso_{S}\mathcal{A}(_{S}S,_{S}S)=S$ such that $f_{i}=\pi _{\overline{G}}$ $%
\alpha _{i}$, and therefore, $Ker$ $f_{i}=$ $Ker$ $\pi _{\overline{H}}\neq
T=Ker$ $\pi _{\overline{H}}$. Thus, $T\neq \cap _{i\in I}Ker$ $f_{i}$.%
\textit{\ \ \ \ \ \ }$_{\square }\medskip \smallskip $

\noindent \textbf{Lemma 3.11.} \textit{Let }$G$\textit{\ and }$H$\textit{\
be subgroups of a group }$S$\textit{, and }$G^{\alpha }=\alpha ^{-1}G\alpha
=H$ \textit{for some }$\alpha \in S$\textit{. Then, }$_{S}\overline{H}\cong $
$_{S}\overline{G}$\textit{\ in }$_{S}\mathcal{A}$\textit{, and therefore,}$\
_{S}\overline{H}\overset{\triangle }{\sim }$\textit{\ }$_{S}\overline{G}$%
\textit{.\medskip }\textbf{\ }

\noindent \textbf{Proof}. Using the equation $G\alpha =\alpha H$, it is
enough only to show that assigning $_{S}\overline{G}$ $\ni sG\longmapsto
s\alpha H\in $ $_{S}\overline{H}$ for any left $G$-coset $sG$ $\in $ $_{S}%
\overline{G}$, one establishes an isomorphism in $_{S}\mathcal{A}$. And this
fact can easily be checked in a habitual fashion.\textit{\ \ \ \ \ \ }$%
_{\square }\medskip \smallskip $

\noindent \textbf{Theorem 3.12.} \textit{For any subgroups }$G$\textit{\ and 
}$H$\textit{\ of a group }$S$\textit{, }$_{S}\overline{H}\overset{\triangle }%
{\sim }$\textit{\ }$_{S}\overline{G}$\textit{\ iff }$G^{\alpha }=H$\textit{\
for some }$\alpha \in S$\textit{.\medskip }\textbf{\ \ }

\noindent \textbf{Proof}.$\Longleftarrow $. This follows from Lemma 3.11.

$\Longrightarrow $. By Propositions 3.9, 3.10 and using Lemma 3.11, one
obtains that $G^{\alpha }=H$.\textit{\ \ \ \ \ \ }$_{\square }\medskip $

\noindent \textbf{Corollary 3.13.} \textit{For normal subgroups }$%
G,H\subseteq S$\textit{,} $\overline{G}\overset{\triangle }{\sim }\overline{H%
}$ \textit{iff }$G=H$\textit{. In particular, for any subgroups }$G,H$%
\textit{\ of an abelian group }$S$\textit{,} $\overline{G}\overset{\triangle 
}{\sim }\overline{H}$ \textit{iff }$G=H$\textit{.\ \ \ \ \ \ }$_{\square
}\medskip $

In what follows, let $z$, perhaps with indexes, denote a singleton
considered as a trivial $S$-act, so to speak, a \textit{zero }$S$-act. Using
this agreement, we obtain\medskip

\noindent \textbf{Proposition 3.14.} \textit{For any }$S$\textit{-act} $A\in
|_{S}\mathcal{A}|$ \textit{without zero subacts, the following statements
are true:}

\textit{(i)} $A\coprod z_{1}\overset{\triangle }{\nsim }A\coprod
(z_{1}\coprod z_{2})$\textit{;}

\textit{(ii) }$A\coprod (z_{1}\coprod z_{2})\overset{\triangle }{\sim }%
A\coprod (\coprod_{i\in I}z_{i})$\textit{, where }$|I|>1$\textit{.} \textbf{%
\ \medskip }

\noindent \textbf{Proof}. (i). Let $F_{X}=S_{x_{1}}\coprod S_{x_{2}}$ for $%
X=\{x_{1},x_{2}\}$, and $\alpha :F_{X}\longrightarrow A\coprod (z_{1}\coprod
z_{2})$ be defined as $\alpha (S_{x_{1}})=z_{1}$ and $\alpha
(S_{x_{2}})=z_{2}$. Then $Ker$ $\alpha =\{(s_{x_{1}},t_{x_{1}})$ $|$ $%
s_{x_{1}},t_{x_{1}}\in S_{x_{1}}\}\cup \{(u_{x_{2}},v_{x_{2}})$ $|$ $%
u_{x_{2}},v_{x_{2}}\in S_{x_{2}}\}\in Cl_{A\amalg (z_{1}\amalg
z_{2})}(F_{X}) $, but it is clear that $Ker$ $\alpha \notin Cl_{A\amalg
z_{1}}(F_{X})$.

(ii). Since $|I|>1$, there is an obvious embedding $A\coprod (z_{1}\coprod
z_{2})\rightarrowtail A\coprod (\coprod_{i\in I}z_{i})$. Consider $S$-act $%
[A\coprod (z_{1}\coprod z_{2})]^{I}$ and the embedding

\noindent $A\coprod (\coprod_{i\in I}z_{i})\rightarrowtail \lbrack A\coprod
(z_{1}\coprod z_{2})]^{I}$ defined as $a\longmapsto (a,a,\ldots ,a,\ldots )$
for any $a\in A$, and let $z_{i}$ go to the string having $z_{1}$ in the $i$%
-th place and $z_{2}$ everywhere else for any $i\in I$. Then, applying
Corollary 3.5, we end the proof.\textit{\ \ \ \ \ \ }$_{\square }\medskip $

Introducing $A^{(\ast )I}\overset{def}{=}\coprod_{i\in I}A_{i},A_{i}=$ $%
_{S}A $, $i\in I$, for any $_{S}A\in |_{S}\mathcal{A}|$ and any set $I$, we
also have\medskip

\noindent \textbf{Proposition 3.15.} \textit{For any }$S$\textit{-act} $A\in
|_{S}\mathcal{A}|$ \textit{without zero subacts, }$B\in |_{S}\mathcal{A}|$%
\textit{, and }$I$ \textit{with} $|I|\geq 1$\textit{, the following
statements are true:}

\textit{(i) }$B$ $\coprod $ $A^{(\ast )I}\coprod z_{1}\overset{\triangle }{%
\sim }B\coprod A\coprod z_{1}$\textit{;}

\textit{(ii) }$B$ $\coprod $ $A^{(\ast )I}\coprod $ $(z_{1}\coprod z_{2})%
\overset{\triangle }{\sim }B\coprod A\coprod $ $(z_{1}\coprod z_{2})$\textit{%
.} \textbf{\ \medskip }

\noindent \textbf{Proof}. (i). Consider $S$-act $[B\coprod A\coprod
z_{1}]^{I}$ and the embedding $B$ $\coprod $ $A^{(\ast )I}\coprod
z_{1}\rightarrowtail \lbrack B\coprod A\coprod z_{1}]^{I}$ such that $a\in $ 
$A_{i}=$ $A\subseteq A^{(\ast )I}$ goes to the string having $a$ in the $i$%
-th place and $z_{1}$ everywhere else for any $i\in I$, and $%
z_{1}\longmapsto (z_{1},z_{1},\ldots ,z_{1},\ldots )$, $b\longmapsto
(b,b,\ldots ,b,\ldots )$ for any $b\in B$. From this, taking into
consideration an obvious embedding $B\coprod A\coprod z_{1}\rightarrowtail B$
$\coprod $ $A^{(\ast )I}\coprod z_{1}$ and applying Corollary 3.5, we get
the statement.

(ii). Just using the embedding $B$ $\coprod $ $A^{(\ast )I}\coprod $ $%
(z_{1}\coprod z_{2})\rightarrowtail \lbrack B\coprod A\coprod $

\noindent $(z_{1}\coprod z_{2})]^{I}$ such that $a\in $ $A_{i}=$ $A\subseteq
A^{(\ast )I}$ goes to the string having $a$ in the $i$-th place and $z_{1}$
everywhere else for any $i\in I$, and $z_{1}\longmapsto (z_{1},z_{1},\ldots
,z_{1},\ldots )$, $z_{2}\longmapsto (z_{2},z_{2},\ldots ,z_{2},\ldots )$,
and $b\longmapsto (b,b,\ldots ,b,\ldots )$ for any $b\in B$, we end the
proof in the similar way as in (i).\textit{\ \ \ \ \ \ }$_{\square }\medskip 
$

\noindent \textbf{Proposition 3.16. }$B$ $\coprod $ $\overline{N}^{(\ast )I}%
\overset{\triangle }{\sim }B$ $\coprod \overline{N}$ \textit{and} $\overline{%
N}^{(\ast )I}\overset{\triangle }{\sim }\overline{N}$ \textit{for any proper
normal subgroup }$N\vartriangleleft S$ \textit{of a group} $S$\textit{, }$%
B\in |_{S}\mathcal{A}|$\textit{, and }$I$ \textit{with} $|I|\geq 1$\textit{%
.\medskip }

\noindent \textbf{Proof}. Since $N$ is a normal subgroup, one can readily
get the embedding $B$ $\coprod $ $\overline{N}^{(\ast )I}\rightarrowtail
\lbrack B\coprod \overline{N}]^{I}$ such that $s\overline{1}\in $ $\overline{%
N}_{i}=$ $\overline{N}\subseteq \overline{N}^{(\ast )I}$ goes to the string
having $s\overline{1}$ in the $i$-th place and $s\overline{t}$ everywhere
else for any $i\in I$, where $\overline{1}$, $\overline{t}\in \overline{N}$
are the cosets modulo $N$ containing $1\in S$ and $t\in S$ , respectively,
and $\overline{1}\neq \overline{t}$, and $b\longmapsto (b,b,\ldots ,b,\ldots
)$ for any $b\in B$. From this, taking into consideration an obvious
embedding $B\coprod \overline{N}\rightarrowtail $ $B$ $\coprod $ $\overline{N%
}^{(\ast )I}$ and applying Corollary 3.5, we get the statement.\textit{\ \ \
\ \ \ }$_{\square }\medskip $

\noindent \textbf{Proposition 3.17.} \textit{For any }$S$\textit{-acts} $%
A,B\in |_{S}\mathcal{A}|$ \textit{without zero subacts, the following
statements are true:}

\textit{(i)} $A\overset{\triangle }{\nsim }$ $B\coprod z$\textit{;}

\textit{(ii)} $A\overset{\triangle }{\nsim }B\coprod $ $(z_{1}\coprod z_{2})$%
\textit{.} \textbf{\ \medskip }

\noindent \textbf{Proof}. (i). Let $F_{X}=S_{x_{1}}\coprod S_{x_{2}}$ for $%
X=\{x_{1},x_{2}\}$, and $\alpha :F_{X}\longrightarrow B\coprod z$ be such a
homomorphism that $\alpha (S_{x_{1}})\subseteq B$ and $\alpha (S_{x_{2}})=z$%
. It is obvious that a non-universal congruence $Ker$ $\alpha \in Cl_{B\text{
}\amalg \text{ }z}(F_{X})$; however, $Ker$ $\alpha \notin Cl_{A}(F_{X})$.

(ii). The same arguments as in (i) work well for this case, too.\textit{\ \
\ \ \ \ }$_{\square }\medskip $

\noindent \textbf{Proposition 3.18.} \textit{For any }$S$\textit{-act} $A\in
|_{S}\mathcal{A}|$ \textit{and }$I$ \textit{with} $|I|\geq 2$\textit{, }$%
A^{(\ast )I}\overset{\triangle }{\sim }A^{(\ast )2}$\textit{, where }$%
A^{(\ast )2}=A_{1}\coprod A_{2}$\textit{\ with }$A_{1}=A=A_{2}$\textit{.} 
\textbf{\ \medskip }

\noindent \textbf{Proof}. Consider the embedding $A^{(\ast
)I}\rightarrowtail (A^{(\ast )2})^{I}=(A_{1}\coprod A_{2})^{I}$ such that $%
a\in $ $A_{i}=$ $A\subseteq A^{(\ast )I}$ goes to the string having $a\in $ $%
A_{1}$ in the $i$-th place and $a\in $ $A_{2}$ everywhere else for any $i\in
I$. Then, using the obvious embeding $A^{(\ast )2}\rightarrowtail $ $%
A^{(\ast )I}$ and applying Corollary 3.5, we get the statement.\textit{\ \ \
\ \ \ }$_{\square }\medskip $

\noindent \textbf{Proposition 3.19. } \textit{Let for }$S$\textit{-acts} $%
A,B\in |_{S}\mathcal{A}|$\textit{\ there exist embedings }$A\overset{i_{A}}{%
\rightarrowtail }B^{I}$\textit{\ and }$B\overset{i_{B}}{\rightarrowtail }%
A^{J}$\textit{\ for some sets }$I$\textit{\ and }$J$\textit{. Then, }$%
C\coprod A\overset{\triangle }{\sim }C\coprod B$\textit{\ for any }$C\in
|_{S}\mathcal{A}|$\textit{.} \textbf{\ \medskip }

\noindent \textbf{Proof}. Consider the following embeddings $C\coprod A%
\overset{1_{C}\text{ }\amalg \text{ }i_{A}}{\rightarrowtail }$\textit{\ }$%
C\coprod B^{I}\overset{\bigtriangleup \text{ }\times \text{ }1}{%
\rightarrowtail }$ $(C\coprod B)^{I}$, where the last embedding defined as
follows: $C\ni c\longmapsto (c,c,\ldots )\in $ $(C\coprod B)^{I}$ and $%
C\coprod B^{I}\ni (b_{i})_{i\in I}\longmapsto (b_{i})_{i\in I}\in (C\coprod
B)^{I}$. Then, considering the similarly defined embeddings $C\coprod B%
\overset{1_{C}\text{ }\amalg \text{ }i_{B}}{\rightarrowtail }$\textit{\ }

\noindent $C\coprod A^{J}\overset{\bigtriangleup \text{ }\times \text{ }1}{%
\rightarrowtail }$ $(C\coprod A)^{J}$ and using Corollary 3.5, we end the
proof.\textit{\ \ \ \ \ \ }$_{\square }\medskip $

\noindent \textbf{Corollary 3.20. }\textit{For any }$S$\textit{-acts} $%
A,C\in |_{S}\mathcal{A}|$ \textit{and }$I$ \textit{with} $|I|\geq 2$\textit{%
, }$C\coprod A^{(\ast )I}\overset{\triangle }{\sim }C\coprod A^{(\ast )2}$%
\textit{.} \textbf{\ \medskip }

\noindent \textbf{Proof}. This follows from the embeddings defined in the
proof of Proposition 3.18 and Proposition 3.19.\textit{\ \ \ \ \ \ }$%
_{\square }\medskip $

The relation $\overset{\triangle }{\sim }$ of \textit{geometric equivalence}
on the class $|_{S}\mathcal{A}|$ of objects of the variety $_{S}\mathcal{A}$
is clearly an equivalence relation. Thus, denoting by \textbf{[}$A$\textbf{]}%
$\overset{def}{=}\{$ $B$ $|B\in |_{S}\mathcal{A}|$ $\&$ $B\overset{\triangle 
}{\sim }A\}$ the class of all $S$-acts geometrically equivalent to an $S$%
-act $A$, we obtain the following a representation theorem for the $\overset{%
\triangle }{\sim }$-equivalence classes of the relation $\overset{\triangle }%
{\sim }$, generalizing Theorem 3.2.\medskip

\noindent \textbf{Theorem 3.21.} \textit{Let }$\{G_{j}\}$ \textit{be} 
\textit{a family of subgroups} $G_{j}\subset S$ \textit{of a group }$S$%
\textit{, }$j\in J$\textit{,\ only containing exactly one \ subgroup\ of
each of the conjugate classes of proper non-normal subgroups of a group }$S$%
\textit{, and }$N_{m}\vartriangleleft S,$\textit{\ }$m\in M$\textit{, all
proper normal subgroups of }$S$\textit{. Then, for any }$S$\textit{-act} $%
A\in $ $|_{S}\mathcal{A}|$\textit{, the }$\overset{\triangle }{\sim }$%
\textit{-equivalence class }\textbf{[}$A$\textbf{]}\textit{\ has a
representation of exactly one of the following three possible types:}

\textit{(i) }\textbf{[}$\coprod_{k\in K}\overline{G_{k}}^{(\ast
)2}\coprod_{t\in T}\overline{G_{t}}\coprod_{l\in L}\overline{N_{l}}$\textbf{]%
}\textit{, where }$K,T$\textit{, and }$L$\textit{\ are some subsets of }$J$%
\textit{\ and }$M$\textit{, respectively, and }$K\cap T=\emptyset $\textit{;}

\textit{(ii) }\textbf{[}$\coprod_{k\in K}\overline{G_{k}}\coprod_{l\in L}%
\overline{N_{l}}\coprod z$\textbf{]}$,K\subseteq J,L\subseteq M$\textit{;}

\textit{(iii) }\textbf{[}$\coprod_{k\in K}\overline{G_{k}}\coprod_{l\in L}%
\overline{N_{l}}\coprod z_{1}\coprod z_{2}$\textbf{]}$,$\textit{\ }$%
K\subseteq J,L\subseteq M$\textit{. \ \medskip }

\noindent \textbf{Proof}. As was mentioned above, from \cite[Propositions
I.5.17 and I.5.34]{klpknmx:maac} one readily gets that any $S$-act $A\in $ $%
|_{S}\mathcal{A}|$ is isomorphic to a coproduct of some $\overline{G_{j}},$ $%
j\in J$, and zero $S$-acts. From this, and using a transfinite induction if
a group $S$ contains an infinite number of subgroups, \textit{i.e.}, $|J\cup
M|\geq \omega $, the result immediately follows from Propositions 3.14,
3.15, 3.16, 3.17, and Corollary 3.20.\textit{\ \ \ \ \ \ }$_{\square
}\medskip $

As a corollary of Theorem 3.21 and Corollary 3.13, we obtain a simpler
version of a representation\ theorem for the $\overset{\triangle }{\sim }$%
-equivalence classes of $S$-acts over abelian groups $S$.\medskip

\noindent \textbf{Theorem 3.22. }\textit{Let }$G_{j}\subset S$\textit{, }$%
j\in J$\textit{, be all proper subgroups of an abelian group }$S$\textit{.
Then, for any }$S$\textit{-act} $A\in $ $|_{S}\mathcal{A}|$\textit{, the }$%
\overset{\triangle }{\sim }$\textit{-equivalence class }\textbf{[}$A$\textbf{%
]}\textit{\ has a representation of exactly one of the following three
possible types:}

\textit{(i) }\textbf{[}$\coprod_{k\in K}\overline{G_{k}}$\textbf{]}\textit{, 
}$K\subseteq J$\textit{;}

\textit{(ii) }\textbf{[}$\coprod_{k\in K}\overline{G_{k}}\coprod z$\textbf{]}%
$,K\subseteq J$\textit{;}

\textit{(iii) }\textbf{[}$\coprod_{k\in K}\overline{G_{k}}\coprod
z_{1}\coprod z_{2}$\textbf{]}$,$\textit{\ }$K\subseteq J$\textit{.\ \ \ \ \
\ }$_{\square }\medskip $

\noindent \textbf{Corollary 3.23.} \textit{In} $_{S}\mathcal{A}$ \textit{%
over an abelian group }$S$\textit{\ of a prime order }$p$\textit{, there are
only the following three }$\overset{\triangle }{\sim }$\textit{-equivalence
classes: }\textbf{[}$S$]\textit{, }[$S\coprod z$\textbf{]}\textit{, and }[$%
S\coprod z_{1}\coprod z_{2}$\textbf{]}\textit{.\ \ \ \ \ \ }$_{\square
}\medskip $

\noindent \textbf{Remark 3.24.} Let $S$ be an abelian group of a prime order 
$p$. Then, $S\coprod z$ is an injective envelope (see, \cite[Section 3.1]%
{klpknmx:maac}) of $S$, however, $S\overset{\triangle }{\nsim }S\coprod z$.
Moreover, from Theorem 3.21 it is easy to see that for any $S$-act $A\in $ $%
|_{S}\mathcal{A}|$ without zero subacts, $A$ and its injective envelope are
not geometrically equivalent.\medskip

Using \cite[Theorem 3.1.8 and Exercise 3.1.9]{klpknmx:maac}, we obtain
another corollary of Theorem 3.21.\medskip

\noindent \textbf{Corollary 3.25.}\textit{\ Let }$G_{j}\subset S$\textit{, }$%
j\in J$\textit{, be all proper subgroups of a group }$S$\textit{. Then, for
any injective }$S$\textit{-act} $A\in $ $|_{S}\mathcal{A}|$\textit{, the }$%
\overset{\triangle }{\sim }$\textit{-equivalence class }\textbf{[}$A$\textbf{%
]}\textit{\ has a representation of exactly one of the following two
possible types:}

\textit{(i) }\textbf{[}$\coprod_{k\in K}\overline{G_{k}}\coprod z$\textbf{]}$%
,K\subseteq J$\textit{;}

\textit{(ii) }\textbf{[}$\coprod_{k\in K}\overline{G_{k}}\coprod
z_{1}\coprod z_{2}$\textbf{]}$,$\textit{\ }$K\subseteq J$\textit{.\ \ \ \ \
\ }$_{\square }\medskip $

Using Proposition 3.16, we also have\medskip

\noindent \textbf{Corollary 3.26.}\textit{\ In} $_{S}\mathcal{A}$\textit{,
all free }$S$\textit{-acts are geometrically equivalent, i.e., they are in }%
\textbf{[}$S$], \textit{provided} \textit{\ }$|S|>1$\textit{.\ \ \ \ \ \ }$%
_{\square }\medskip $

In conclusion, we state the following, in our view, interesting and
promising problems.\medskip

\noindent \textbf{Problem 1.} Our conjecture is that in general the \textit{%
\ }$\overset{\triangle }{\sim }$-equivalence\textit{\ }classes\ $[A]$ in
Theorems 3.21, 3.22, and Corollary 3.25 may have different representations
of the same type. Therefore, it would be interesting and important to
describe groups $S$ over which any \textit{\ }$\overset{\triangle }{\sim }$%
-equivalence\textit{\ }class $[A]$, $A\in $ $|_{S}\mathcal{A}|$, has a
unique representation of the given types.\medskip

\noindent \textbf{Problem 2.} In light of Corollaries 3.25 and 3.26, in
which were considered the classes of injective and free $S$-acts,
respectively, it is interesting to obtain analogs of those results for other
homological classes of $S$-acts (see, \textit{e.g.}, \cite{klpknmx:maac}%
).\medskip

\noindent \textbf{Problem 3.} To extend the considerations results of the
present paper, obtained for $S$-acts over groups $S$, to some other
interesting varieties $_{S}\mathcal{A}$ over ``nice'' monoids $S$.\


\begin{thebibliography}{99}
\bibitem{baummysrem:agog} G. Baumslag, A. Myasnikov, and V. N.
Remeslennikov, Algebraic geometry over groups, \textit{J. Algebra} \textbf{%
219} (1999) 16--79.

\bibitem{berzinsplts:agivoaw} A. Berzins, B. Plotkin, and E. Plotkin,
Algebraic geometry in varieties of algebras with the given algebra of
constants, \textit{J. Math. Sci.} \textbf{102} (2000) 4039--4070.

\bibitem{berzins:geoa} A. Berzins, Geometric equivalence of algebras, 
\textit{Internat. J. Algebra Comput.} \textbf{11} (2001) 447--456.

\bibitem{gobsh:rappcgeg} R\"{u}diger G\"{o}bel and Saharon Shelah, Radicals
and Plotkin's problem concerning geometrically equivalent groups, \textit{%
Proc. Amer. Math. Soc.} \textbf{130} (2002) 673--674.

\bibitem{gratzr:ua} G. Gr\"{a}tzer, \textit{Universal Algebra}, 2nd Ed.
(Springer-Verlag, New York-Berlin, 1979).

\bibitem{katlipplot:aocofmfsflm} Y. Katsov, R. Lipyanski, and B. Plotkin,
Automorphisms of categories of free modules, free semimodules, and free Lie
modules, \textit{Comm. Algebra} \textbf{35} (2007) 931--952.

\bibitem{kat:apop} Y. Katsov, A problem of B. Plotkin for $S$-acts:
automorphisms of categories of free $S$-acts, \textit{Comm. Algebra} \textbf{%
35} (2007) 1709--1714.

\bibitem{klpknmx:maac} M. Kilp, U. Knauer, and A. V. Mikhalev, \textit{%
Monoids, Acts and Categories} (Walter de Gruyter, Berlin-New York, 2000).

\bibitem{macl:cwm} S. Mac Lane, \textit{Categories for the Working
Mathematician} (Springer-Verlag, New York-Berlin, 1971).

\bibitem{mal:as} A. I. Maltsev, \textit{Algebraic Systems} (Springer-Verlag,
New York-Berlin, 1973).

\bibitem{plotkin:voaaavcoav} B. Plotkin, Varieties of algebras and algebraic
varieties. Categories of algebraic varieties, \textit{Siberian Adv. Math.} 
\textbf{7} (1997) 64--97.

\bibitem{plotkin:snoagiua} B. I. Plotkin, Some notions of algebraic geometry
in universal algebra, (Russian) \textit{Algebra i Analiz }\textbf{9} (1997)
224--248; English transl., \textit{St. Petersburg Math.} J.\textbf{\ 9}
(1998) 859--879.

\bibitem{plplts:geog} B. Plotkin, E. Plotkin, A. Tsurkov, Geometrical
equivalence of groups, \textit{Comm. Algebra} \textbf{27} (1999) 4015--4025.

\bibitem{plotkin:slotuag} B. Plotkin, Seven lectures on the universal
algebraic geometry, \textit{Preprint},\textit{\ Institute of Mathematics,
Hebrew University, Jerusalem}, arXiv:math. GM/0204245 (2002).

\bibitem{plotkin:piaibuag} B. I. Plotkin, Problems in algebra inspired by
universal algebraic geometry, (Russian) \textit{Fundam. Prikl. Mat.} \textbf{%
10} (2004) 181--197.

\bibitem{plts:atgeorog} B. Plotkin, A.Tsurkov, Action type geometrical
equivalence of representations of groups, \textit{Preprint} arXiv:math.
RT/0501337 (2005).

\bibitem{tsurkov:geontfg} A. Tsurkov, Geometrical equivalence of nilpotent
torsion free groups, \textit{Preprint} arXiv:math. GR/0411313 (2004).
\end{thebibliography}
\end{document}